\documentclass[12pt]{amsart}

\usepackage[dvips]{graphicx}
\usepackage{amsmath,amssymb,amsthm,wrapfig,amsfonts,enumerate,latexsym}

\usepackage{comment,color}

\newtheorem{theorem}{Theorem}[section]
\newtheorem{lemma}[theorem]{Lemma}
\newtheorem{corollary}[theorem]{Corollary}
\newtheorem{proposition}[theorem]{Proposition}

\theoremstyle{definition}

\newtheorem{definition}[theorem]{Definition}
\newtheorem{remark}[theorem]{Remark}
\newtheorem{condition*}[theorem]{Condition}
\newtheorem*{acknowledgement}{Acknowledgements}       
 
\numberwithin{equation}{section}

\begin{document}

\title[Fixed point property for a ${\rm CAT}(0)$ space]{Fixed point property for a {\boldmath $\mathrm{CAT}(0)$} space 
which admits a proper cocompact group action}
\author
[T. Toyoda]{Tetsu Toyoda}

\email
        [Tetsu Toyoda]{toyoda@genl.suzuka-ct.ac.jp}

\address
        [Tetsu Toyoda]
        {\endgraf Suzuka National College of Technology, 
        \endgraf Shiroko-cho, Suzuka, Mie, 510-0294, Japan}

\date{}


\begin{abstract}
We prove that if a geodesically complete 
$\mathrm{CAT}(0)$ space $X$ admits a proper cocompact isometric action of a group, 
then  the Izeki-Nayatani invariant of $X$ is less than $1$. 
Let $G$ be a finite connected graph, $\mu_1 (G)$ be the linear spectral gap of $G$, and $\lambda_1 (G,X)$ 
be the nonlinear spectral gap of $G$ with respect to such a $\mathrm{CAT}(0)$ space $X$. 
Then, the result implies that the ratio $\lambda_1 (G,X) / \mu_1 (G)$ 
is bounded from below by a positive constant which is independent of 
the graph $G$. 
It follows that any isometric action of a random group of the graph model on such $X$ has a global 
fixed point. 
In particular, any isometric action of a random group of the graph model on 
a Bruhat-Tits building associated to a semi-simple algebraic group has a global fixed point. 
\end{abstract}

\maketitle

\section{Introduction}\label{intro-sec}

\subsection{Nonlinear spectral gaps}

Let $G=(V,E)$ be a graph, 
where $V$ and $E$ denote the sets 
of vertices and unoriented edges, respectively. 
Throughout this paper, 
we assume that every graph is simple and connected and  
satisfies $2\leq|V|<\infty$. 
A {\em weight function} on $G$ is a symmetric function 
$m :V\times V\to\lbrack 0,\infty )$ 
whose support equals the set $\vec{E}=\{ (u,v)\in V\times V\hspace{1mm}|\hspace{1mm}\{ u,v\}\in E\}$. 
A weight function $m$ induces a weight 
$m(u)$ of each vertex $u\in V$ by $m(u)=\sum_{v\in V}m(u,v)$. 
We use the convention that $m(V )=\sum_{v\in V}m(v)$. 
The pair $(G,m)$ is called a \textit{weighted graph}. 
Unless we specify otherwise, we assume that every graph  
is equipped with the {\em uniform} weight function $m$ defined as 
$$
m(u,v)=
\begin{cases}
1,&\quad \textrm{if}\quad (u,v)\in \vec{E}, \\
0,&\quad \textrm{otherwise} .
\end{cases}
$$

The {\em linear spectral gap} $\mu_1 (G)$ of a weighted graph $(G,m)$ is the 
first positive eigenvalue of the combinatorial Laplacian $\Delta$ which acts on 
functions $f:V\to\mathbb{R}$ as 
$$
\Delta f (v)
=
f(v)-
\sum_{u\in V}\frac{m(v,u)}{m(v)}f(u), \quad 
v\in V. 
$$
It can be computed variationally as 
\begin{equation}\label{RQ-eq}
\mu_1 (G)=\inf\left\{\frac{\frac{1}{2}\sum_{u,v\in V}m(u,v)|f(u)-f(v)|^2}{\sum_{v\in V}m(v)|f(v)-\overline{f}|^2}\hspace{1mm}\bigg|\hspace{1mm}
f:V\to\mathbb{R}\textrm{ is nonconstant}\right\} , 
\end{equation}
where $\overline{f}=\{1 /m(V)\}\sum_{v\in V}\{m(x)f(x)\}$, or 
\begin{equation}\label{RQgro-eq}
\mu_1 (G)=\inf\left\{\frac{\sum_{u,v\in V}m(u,v)|f(u)-f(v)|^2}{\frac{1}{m(V)}\sum_{u,v\in V}m(u)m(v)|f(u)-f(v)|^2}\hspace{1mm}\bigg|\hspace{1mm}
f:V\to\mathbb{R}\textrm{ is nonconstant}\right\} . 
\end{equation}

Recently, several nonlinear analogues of $\mu_1 (G)$ with respect to a general metric space $X$ were defined 
by considering mappings $f:V\to X$ instead of $\mathbb{R}$-valued functions. 
They are called {\em nonlinear spectral gaps} and 
played important roles in metric geometry and geometric group theory.

By generalizing the formula \eqref{RQ-eq}, we obtain the following definition of a nonlinear spectral gap, which was first 
introduced by M. -T. Wang \cite{Wa} for the case where the target metric space is an Hadamard manifold.  
Throughout this paper, every metric space is assumed to contain at least two points. 

\begin{definition}[Wang invariant]\label{wang-def}
Let $\left(G,m\right)$ be a weighted graph and $(X,d_X)$ be 
a complete $\mathrm{CAT}(0)$ space. 
The {\em Wang invariant} $\lambda_1 (G,X)$ of $G$ with respect to $X$ is defined as 
\begin{equation*}
\lambda_1 (G,X)=\inf\left\{\frac{\frac{1}{2}\sum_{u,v\in V}m(u,v)d_X (f(u),f(v))^2}{\sum_{v\in V}m(v)d_X (f(v),\overline{f})^2}\hspace{1mm}\bigg|\hspace{1mm}
f:V\to X\textrm{ is nonconstant}\right\} , 
\end{equation*}
where $\overline{f}$ denotes the barycenter of the probability measure 
$$
\sum_{v\in V}\frac{m(v)}{m(V )}\mathrm{Dirac}_{f(v)}
$$
on $X$. 
Here, $\mathrm{Dirac}_{f(v)}$ denotes the Dirac measure at $f(v)\in X$. 
\end{definition}

The Wang invariant plays a crucial role in the theory of rigidity of groups (\cite{Wa}, \cite{IN}, \cite{P}). 
By generalizing the formula \eqref{RQgro-eq}, 
we obtain the definition of another nonlinear spectral gap which was defined by Gromov \cite{G}.

\begin{definition}
Let $\left(G,m\right)$ be a weighted graph and $(X,d_X)$ be a metric space. 
The {\em Gromov nonlinear spectral gap} $\lambda_1^{\mathrm{Gro}} (G,X)$ is defined as 
\begin{multline*}
\lambda_1^{\mathrm{Gro}} (G,X) \\
=\inf\left\{\frac{\sum_{u,v\in V}m(u,v)d_X (f(u), f(v))^2}{\frac{1}{m(V)}\sum_{u,v\in V}m(u)m(v)d_X (f(u), f(v))^2}\hspace{1mm}\bigg|\hspace{1mm}
f:V\to X\textrm{ is nonconstant}\right\} .
\end{multline*}
\end{definition}
By definition, we have 
$$
\lambda_1 (G,\mathbb{R})=\lambda_1^{\mathrm{Gro}}(G,\mathbb{R})=\mu_1 (G)
$$
for any graph $G$. 
Moreover, we also have   
$$
\lambda_1 (G,\mathcal{H}) =\lambda_1^{\mathrm{Gro}} (G,\mathcal{H}) =\mu_1 (G). 
$$
for any graph $G$ and Hilbert space $\mathcal{H}$. 
For a general complete $\mathrm{CAT}(0)$ space $X$ and a graph $G$, 
these two nonlinear spectral gaps have the following relation (see \cite{P}): 
\begin{equation}\label{wang-gromov-ineq}
\frac{1}{2}\lambda_1 (G,X)\leq\lambda_1^{\mathrm{Gro}} (G,X)\leq\lambda_1 (G,X) .
\end{equation}

\subsection{Comparison of nonlinear and linear spectral gaps}

It is a fundamental question to ask for what kind of complete $\mathrm{CAT}(0)$ space $X$, 
does there exist a constant $C_X >0$ depending only on $X$ which satisfies 
\begin{equation}\label{ratio-bound}
\lambda_1 (G,X)\geq C_X \mu_1 (G)
\end{equation}
for any graph $G$. 
It is known that the existence of such a constant $C_X >0$ implies many important conclusions 
including the following (A) and (B): 
\begin{enumerate}
\item[(A)]
Any isometric action of a random group of the graph model on $X$ has 
a global fixed point. 
\item[(B)]
A sequence of expanders does not embed coarsely into $X$ (\cite{G}). 
\end{enumerate}
The conclusion (A) was proved by Izeki, Kondo and Nayatani in \cite{IKN2}. 
For its precise statement, see Theorem \ref{IKN-fpp} in Section \ref{applications-sec}. 
Roughly, it states that if we equip a suitable probability measure with 
a set $\mathcal{G}$ of finitely generated groups, then, with high probability, 
a randomly chosen group $\Gamma\in\mathcal{G}$ is infinite and 
any isometric action of $\Gamma$ on $X$ has a global fixed point. 
This guarantees the existence of infinite groups $\Gamma$ whose isometric actions on $X$ always have global fixed points.

For the definitions of coarse embeddings and sequences of expanders, see Section \ref{applications-sec}. 
In \cite{G}, Gromov proved that a sequence of expanders does not embed coarsely into a Hilbert space. 
Since then, coarse embeddability of a sequence of expanders into a metric space has become 
an important obstruction of the space to be embedded coarsely into a Hilbert space. 
The conclusion (B) states that $X$ does not have such an obstruction to embed coarsely into a Hilbert space, 
and is proved easily by applying the argument of Gromov. 
For the detailed proof of the conclusion (B), see Theorem 4.5 of \cite{FT}.

The purpose of this paper is to specify complete $\mathrm{CAT}(0)$ spaces $X$ 
which allow the existence of such constants $C_X >0$ as above. 
Consequently, it will specify spaces which satisfy the above conclusions (A) and (B). 

Throughout this paper, we denote by $B(p,r)$ the open ball of radius $r$ centered at $p$ and 
by $\overline{B}(p,r)$ the closed ball of radius $r$ centered at $p$.  
We use the following definition. 
\begin{definition}[\cite{BH}, Chapter I.8]\label{cocompact-def}
An isometric action of a group $\Gamma$ on a metric space $X$ is called \textit{cocompact} if 
there exists a compact subset $K\subset X$ such that $X=\cup_{\gamma\in\Gamma}\gamma K$. 
An isometric action of $\Gamma$ on a metric space $X$ 
is called \textit{proper} if for each $p\in X$ there exists $r>0$ such that 
the set 
$\{\gamma\in\Gamma\hspace{1mm}|\hspace{1mm}\gamma B(p,r)\cap B(p,r)\neq\phi\}$ 
is finite. 
\end{definition}
We prove the following theorem. 
\begin{theorem}\label{cocompact-th}
Let $\mathcal{X}=\{ X_1 ,X_2 ,\ldots ,X_n \}$ be a finite set of 
geodesically complete $\mathrm{CAT}(0)$ spaces such that 
each $X_i$ admits a proper cocompact isometric action of a group. 
Then, there exists a constant $C=C_{\mathcal{X}} >0$ which depends only on $\mathcal{X}$ such that the inequality 
\begin{equation}\label{linear-bound-ineq-cocompact} 
\lambda_1 (G, X)\geq C \mu_1 (G) 
\end{equation}
holds whenever $X$ is a (finite or infinite) product of 
copies of spaces in $\mathcal{X}$ and $G$ is a graph. 
In particular, any Bruhat-Tits building $X$ associated to a semi-simple algebraic group 
admits the existence of a constant $C=C_{X}$ which satisfies the inequality \eqref{linear-bound-ineq-cocompact} for any graph $G$. 
\end{theorem}
We also prove an estimate of the same type when a complete $\mathrm{CAT}(0)$ space $X$  
is uniformly locally doubling in the following sense. 

\begin{definition}\label{doubling-def}
Fix $N\in\lbrack 1,\infty )$. 
A metric space is 
called {\em doubling with doubling constant} $N$ 
if every closed ball  
can be covered by at most $N$ closed balls of half the radius. 
We say that a metric space is {\em uniformly locally doubling with doubling constant} 
$N$ 
if any point has a neighborhood which is doubling with doubling constant $N$. 
\end{definition}
\begin{theorem}\label{doubling-th}
For each $N\in\lbrack 1,\infty )$,  
there exists a constant $C=C_N >0$ such that the inequality 
\begin{equation}\label{linear-bound-doubling-ineq} 
\lambda_1 (G, X)\geq C \mu_1 (G)
\end{equation}
holds for every graph $G$ and a complete $\mathrm{CAT}(0)$ space $X$ which is isometric to 
a (finite or infinite) product of uniformly locally doubling $\mathrm{CAT}(0)$ spaces with a 
common doubling constant $N$. 
\end{theorem}
Theorem \ref{cocompact-th} and Theorem \ref{doubling-th} yield that 
if a complete $\mathrm{CAT}(0)$ space $X$ satisfies the hypothesis of either theorem, 
$X$ satisfies the conclusions (A) and (B). 
We state this explicitly as Threorem \ref{fpp-conclusion-th} and Theorem \ref{coarse-conclusion-th} in Section \ref{applications-sec}. 
In particular, if $X$ is a Bruhat-Tits building associated to a semi-simple algebraic group, 
then any isometric action of a random group of the graph model on $X$ 
has a global fixed point.

\subsection{Relations with other results}

Naor-Silberman \cite{NS} proved that 
if a metric space $X$ has finite Nagata dimension, then 
for every $\varepsilon >0$, 
there exists a constant $C_{X,\varepsilon}$ which satisfies 
\begin{equation}\label{NS-estimate}
\lambda_1^{\mathrm{Gro}}(G ,X)\geq C_{X,\varepsilon} \mu_1 (G)^{1+\varepsilon} ,
\end{equation}
for every graph $G$. 
Moreover, 
Naor-Silberman \cite{NS} also proved that 
the weaker inequality \eqref{NS-estimate} suffices to imply the 
fixed point property (A) of a random group for $X$ whenever $X$ is $p$-uniformly convex for some $p\geq 2$. 
Since each Bruhat-Tits building associated to a semi-simple algebraic group has finite Nagata dimension, 
and complete $\mathrm{CAT}(0)$ spaces are $2$-uniformly convex, 
for such a building $X$, 
the fixed point property (A) also follows from their result.

However, an advantage of our result is that our estimate \eqref{ratio-bound} is better than their estimate \eqref{NS-estimate}. 
In fact, by \eqref{wang-gromov-ineq}, we obtain the following corollaries of Theorem \ref{cocompact-th} and Theorem \ref{doubling-th}, respectively, 
which improve Naor-Silbermann's estimate \eqref{NS-estimate} 
when the target metric space satisfies the hypothesis of either theorem.  

\begin{corollary}
If a 
geodesically complete 
$\mathrm{CAT}(0)$ space $X$ admits a proper cocompact isometric action of a group, 
then there exists a constant $C'_X >0$ such that the inequality 
\begin{equation}\label{linear-bound-ineq} 
\lambda_1^{\mathrm{Gro}} (G, X)\geq C'_X \mu_1 (G) 
\end{equation}
holds for every graph $G$. 
\end{corollary}

\begin{corollary}
For each $N\in\lbrack 1,\infty )$,  
there exists a constant $C'_N >0$ such that the inequality 
\begin{equation}\label{linear-bound-doubling-ineq} 
\lambda_1^{\mathrm{Gro}} (G, X)\geq C'_N \mu_1 (G)
\end{equation}
holds for every uniformly locally doubling complete $\mathrm{CAT}(0)$ space $X$ 
with doubling constant $N$ and every graph $G$.
\end{corollary}

\subsection{The Izeki-Nayatani invariant}
To obtain such a stronger estimate, we use the so-called Izeki-Nayatani invariant. 
Izeki and Nayatani introduced the Izeki-Nayatani invariant $0\leq\delta (X)\leq 1$ of 
a complete $\mathrm{CAT}(0)$ space $X$ in \cite{IN}, 
and proved that 
\begin{equation}\label{delta-wang}
\lambda_1 (G,X)\geq\left(1-\delta (X)\right)\mu_1 (G)
\end{equation}
for any complete $\mathrm{CAT}(0)$ space $X$ and graph $G$. 
For the definition of the Izeki-Nayatani invariant, see Section \ref{delta-sec}. 

On the other hand, a standard method to compare the linear spectral gap and the nonlinear spectral gap with respect to $X$  
is to estimate the bi-Lipschitz distortion of $X$ into a Hilbert space. 
The bi-Lipschitz distortion $c_2 (X)$ of a metric space $X$ into a Hilbert space is 
the infimum of $D>0$ such that there exists a $1$-Lipschitz mapping $f:X\to\mathcal{H}$ to a Hilbert space which satisfies 
$$
\frac{1}{D}d_X (x,y)\leq\| f(x)-f(y)\|\leq d_X (x,y)
$$
for every $x,y\in X$, and we have 
$$
\lambda_1^{\mathrm{Gro}} (G,X)\geq\frac{1}{c_2 (X)^2} \mu_1 (G). 
$$

Although the Izeki-Natayatani invariant $\delta (X)$ can also be estimated by using the bi-Lipschitz distortion $c_2 (X)$ into a Hilbert space,  
the present author \cite{To2} established another method to estimate it, which even does not  
require the existence of bi-Lipschitz embeddings of $X$ into a Hilbert space. 
This method enables us to obtain our estimates. 
We summarize this method in Section \ref{delta-sec}.

\subsection{Organization}
The paper is organized as follows. 
In Section \ref{preliminaries-sec}, we briefly review some basic notions concerning $\mathrm{CAT}(0)$ spaces. 
In Section \ref{delta-sec}, we recall the definition of the Izeki-Nayatani invariant and discuss some basic properties of it. 
We also summarize the method obtained in \cite{To2} to estimate this invariant. 
In Section \ref{cocompact-sec}, we prove Theorem \ref{cocompact-th}. 
In Section \ref{doubling-sec}, we prove Theorem \ref{doubling-th}. 
To prove Theorem \ref{doubling-th}, we prove that 
the ultralimit of a sequence of doubling length spaces with a common doubling constant is also 
doubling with the same constant. 
In Section \ref{applications-sec}, we see that our results imply fixed-point theorems of random groups and 
non-embeddability of sequences of expanders. 
In Appendix, we discuss some other facts concerning the Izeki-Nayatani invariant.

\begin{acknowledgement}
I would like to thank S. Nayatani, K. Fujiwara, H. Izeki, and T. Kondo 
for helpful discussions. 
\end{acknowledgement}

\section{Preliminaries}
\label{preliminaries-sec}
In this section, we briefly recall some basic notions in metric geometry. 
For a detailed exposition, 
we refer the reader to \cite{BH}, \cite{BBI}, and \cite{Ro}. 

Let $(X,d_X)$ be a metric space. 
A continuous mapping $\gamma :I\to X$ from an interval $I\subset\mathbb{R}$ to $X$ is called a {\em path} in $X$. 
When $I=\lbrack a,b\rbrack$ is a closed interval, it is called a path joining $\gamma (a)$ to $\gamma (b)$. 
The {\em length} $L(\gamma )$ of a path $\gamma :\lbrack a,b\rbrack\to X$ is defined as 
$$
L(\gamma )=\sup\sum_{i=1}^k d_X \left(\gamma (t_{i-1}),\gamma(t_i) \right), 
$$
where the supremum is taken over all finite subdivisions 
$$
a= t_0 \le t_1 \le\cdots\le t_k=b .
$$
A path $\gamma :\lbrack a,b\rbrack\to X$ is called {\em arc-length parametrized} if 
$L(\gamma |_{\lbrack a,t\rbrack})=|t-a|$ for all $t\in\lbrack a,b\rbrack$, where 
$\gamma |_{\lbrack a,t\rbrack}$ is the restriction of $\gamma$ to $\lbrack a,t\rbrack$. 
Any path can be reparametrized to an arc-length parametrized path. 
$X$ is called a {\em length space} if the distance 
$d_X (p,q)$ between any two points $p,q\in X$ 
is equal to the infimum over the lengths of paths joining $p$ to $q$. 
We call a path $\gamma :I\to X$ a {\em geodesic} 
if it is an isometric embedding of the interval $I$ to $X$. 
A metric space is called a {\em geodesic space} if every pair of points is joined by a geodesic. 
We call a path $\gamma :I\to X$ a {\em local geodesic} 
if for every $t\in I$ there exists a neighborhood $J$ of $t$ in $I$ such that 
the restriction $\gamma |_{J} :J\to Y$ is 
a geodesic. 
\begin{definition}
A metric space $X$ is called {\em geodesically complete} if 
it is complete and 
any local geodesic $\gamma :\lbrack 0,a\rbrack\to X$ is a restriction of some local geodesic 
$\tilde{\gamma}:\lbrack 0,b\rbrack\to X$ with $0<a<b$. 
\end{definition}

A {\em geodesic triangle} in $X$ is a triple 
$\triangle =(\gamma_1 ,\gamma_2 ,\gamma_3 )$ 
of geodesics 
$\gamma_i :\lbrack a_i,b_i \rbrack\to X$ 
such that 
$$
\gamma_1 (b_1)=\gamma_2 (a_2),\quad\gamma_2 (b_2)=\gamma_3 (a_3), 
\quad\gamma_3 (b_3)=\gamma_1 (a_1) .
$$
For a geodesic triangle $\triangle =(\gamma_1 ,\gamma_2 ,\gamma_3 )$ 
there is a geodesic triangle 
$$
\overline{\triangle}=
(\overline{\gamma_1} ,\overline{\gamma_2} ,
\overline{\gamma_3} ),\quad
\overline{\gamma_i} :\lbrack a_i,b_i \rbrack\to\mathbb{R}^2
$$
in $\mathbb{R}^2$
such that $L(\gamma_i )=L(\overline{\gamma_i})$ for each $i$. 
Such a triangle $\overline{\triangle}$ is unique up to isometry of $\mathbb{R}^2$.  
We call it the {\em comparison triangle} of $\triangle$ in 
$\mathbb{R}^2$. 
A geodesic triangle $\triangle$ is said to be \textit{thin} 
if 
$$
d_Y (\gamma_i (s),\gamma_j (t))\leq
d_{\kappa}(\overline{\gamma_i} (s),\overline{\gamma_j} (t))
$$
whenever $i,j\in\{ 1,2,3\}$, $s\in\lbrack a_i ,b_i \rbrack$, 
and $t\in\lbrack a_j ,b_j \rbrack$. 

\begin{definition}
A geodesic space $X$ is called a $\mathrm{CAT}(0)$ space 
if every geodesic triangle in $X$ is thin.  
\end{definition}

By definition, for any pair of points $p,q\in X$, 
a geodesic $\gamma :\lbrack 0,d_X (p,q)\rbrack\to X$ joining $p$ to $q$ 
is unique whenever $X$ is a $\mathrm{CAT}(0)$ space. 
It is known that 
every local geodesic $\gamma :\lbrack a,b \rbrack\to X$ in a geodesically complete space  
is a restriction of 
a local geodesic 
$\tilde{\gamma} :\mathbb{R}\to Y$ (see \cite[Corollary 9.1.28.]{BBI}), 
and every local geodesic in a complete $\mathrm{CAT}(0)$ space is a geodesic (see \cite[Chapter II, Proposition 1.4]{BH}). 
Thus, every geodesic 
$\gamma :\lbrack a,b\rbrack\to X$ in a geodesically complete $\mathrm{CAT}(0)$ space $X$ 
is a restriction of a geodesic $\tilde{\gamma}:\mathbb{R}\to X$.

Let $\gamma :\lbrack a,b\rbrack\to X$, 
$\gamma' :\lbrack a',b' \rbrack\to X$
be two geodesics in a $\mathrm{CAT}(0)$ space $X$ with  
$\gamma (a)=\gamma' (a')=p$. 
We define the \textit{angle} $\angle_p (\gamma ,\gamma' )$ 
between $\gamma$ and $\gamma'$ as 
$$
\angle_p (\gamma ,\gamma' )=
\lim_{t\to a,t'\to a'}\angle^{0}_p (\gamma(t),\gamma' (t') ) ,
$$ 
where $\angle^{0}_p (\gamma(t),\gamma'(t') ) $ is 
the corresponding angle of the triangle in $\mathbb{R}^2$ 
whose side lengths are $d_X( p, \gamma (t))$, $d_X (\gamma (t),\gamma' (t))$ and $d_X (\gamma' (t), p)$. 
The existence of the limit is guaranteed by the definition of 
$\mathrm{CAT}(0)$ spaces. 
The law of cosines on a Euclidean space yields 
\begin{equation}\label{cosine}
\cos \angle_p (\gamma ,\gamma' )
=
\lim_{t\to a, t'\to a'}\frac{d_X (p,\gamma (t))^2 + d_X (p,\gamma' (t') )^2 
-d_X (\gamma (t) ,\gamma' (t'))^2}{2 d_X (p,\gamma (t)) d_X (p,\gamma' (t') )}. 
\end{equation}

\begin{definition}\label{euclidean-cone}
Let $(S,d_S )$ be a metric space. 
The cone $\mathrm{Cone}(S)$ over $S$ is 
the quotient of the product $S\times\lbrack 0,\infty )$ 
obtained by identifying all points in 
$S\times\{ 0\}\subset S\times\lbrack 0,\infty )$. 
The point represented by $(x,0)$ for any $x\in S$ is called the \textit{origin} of the cone and 
we denote this point by $o$.  
The cone distance $d_{\mathrm{Cone}(S)} (v,w)$ 
between two points $v,w\in\mathrm{Cone}(S)$ 
represented by $(x,t),(y,s)\in S\times\lbrack 0,\infty )$ 
respectively, is defined by 
$$
d_{\mathrm{Cone}(S)} (v,w)= 
\sqrt{t^2 +s^2 -2ts\cos (\min\{\pi ,d_S (x,y)\})} .
$$
Then $(\mathrm{Cone}(S),d_{\mathrm{Cone}(S)})$ becomes a metric space. 
We call this metric space 
the \textit{Euclidean cone} over $(S,d_S )$. 
\end{definition}

For an element $v\in\mathrm{Cone}(S)$ represented by $(x,r)\in S\times\lbrack 0,\infty )$ and $c>0$, 
we denote  by $cv$ the element represented by $(x, cr)$. 
We claim that 
$$
d_{\mathrm{Cone}(S)} (cv,cw)= cd_{\mathrm{Cone}(S)} (v,w)
$$
holds for any $v,w\in\mathrm{Cone}(S)$

\begin{definition}\label{tangent-cone}
Let $(X,d_X )$ be a  
$\mathrm{CAT}(0)$ space, and let $p\in X$. 
We denote by $(S_p X )^{\circ}$ the quotient set of all nontrivial geodesics starting from $p$ by 
the equivalence relation $\sim$ defined by 
$\gamma\sim\gamma'\Longleftrightarrow\angle_p (\gamma ,\gamma')=0$.   
Then the angle $\angle_p$ induces a distance on  
$(S_p X )^{\circ}$, which we denote by the same symbol $\angle_p$. 
The {\em space of directions} $S_p X$ at $p$ is 
the metric completion of the metric space $\left( (S_p X)^{\circ}, \angle_p \right)$. 
The {\em tangent cone} $TC_p X$ of $X$ at $p$ is the 
Euclidean cone $\mathrm{Cone}(S_p X )$ over the 
space of directions at $p$. 
Define a map $\pi_p :X\to TC_p X$ by 
$\pi_p (q)=\left(\lbrack \gamma\rbrack ,d_X(p,q) \right)$ 
where $\lbrack\gamma\rbrack$ is the 
equivalence class represented by the unique geodesic $\gamma$ joining $p$ and $q$. 
\end{definition}
It is easily seen that the map $\pi_p$ defined as above is $1$-Lipschitz. 
It is also seen that each tangent cone $TC_p X$ is the metric completion of 
the Euclidean cone $\mathrm{Cone}((S_p X)^{\circ} )$. 
If we denote the canonical inclusion of $S_p X$ into $TC_p X$ by $\iota$, 
then it is straightforward from the definition of the metric on Euclidean cones that 
we have 
\begin{equation}\label{inclusion}
\frac{2}{\pi}d_S(x,y)
\le d_T (\iota (x),\iota (y))
\le d_S(x,y)
\end{equation}
for all $x,y \in S_p X$, where $d_S$ and $d_T$ represent the distance functions of 
$S_p X$ and $TC_p X$, respectively. 

The $\mathrm{CAT}(0)$ condition is preserved under taking ($\ell_2$-)product. 
\begin{definition}
Let 
$(X_1 ,d_1 ) ,(X_2 ,d_2 ), \ldots$ be metric spaces 
with basepoints $o_1 \in X_1$, $o_2 \in X_2$, $\ldots$, respectively. 
The ($\ell^2$-)\textit{product} $X$ of $X_1 , X_2 , \ldots$ with respect to the basepoints $o_1 ,o_2 ,\ldots$
consists of all sequences $(x_n )_n$ with $x_n \in X_n$, satisfying 
$\sum_{n}d_n (o_n ,x_n )^2 <\infty$, 
and is equipped with the metric function $d$ 
defined by 
$$
d(x,y )^2 =\sum_{n=1}^{\infty}d_n (x_n , y_n )^2
$$
for any elements $x=(x_1 ,x_2 ,\ldots )\in X$ and $y=(y_1 ,y_2 ,\ldots )\in X$. 
\end{definition}

To define the Wang invariant, we need to consider a finitely supported probability measure on a complete $\mathrm{CAT}(0)$ space. 
We often write a finitely supported probability measure $\mu$ on a metric space $X$ in the form 
$$
\mu=\sum_{i=1}^m t_i \mathrm{Dirac}_{p_i}, 
$$
where $\mathrm{Dirac}_{p_i}$ is the Dirac measure at $p_i \in X$ and each $t_i$ is 
the weight $\mu (\{p_i \})$ at $p_i$. 
We denote the support of a measure $\mu$ by $\mathrm{Supp}(\mu )$. 
When $X$ is a complete $\mathrm{CAT}(0)$ space, 
there exists a unique point on $X$ which minimizes the function 
$p\mapsto\int_X d_X (p,q)^2 \mu (dq)=\sum_{i=1}^m t_i d_Y (p, p_i )^2 $ 
(see \cite{St}). 
This point is called a {\em barycenter} of $\mu$ and denoted by $\mathrm{bar}(\mu )$.

\section{Izeki-Nayatani invariant}\label{delta-sec}
In this section, we recall 
the definition of the Izeki-Nayatani invariant $\delta$ and 
its basic properties.

\begin{definition}[Izeki-Nayatani \cite{IN}]\label{delta-def}
Let $X$ be a complete $\mathrm{CAT} (0)$ space, 
and $\mathcal{P} (X)$ be the space of all finitely supported probability measures $\mu$ with 
$|\mathrm{Supp}(\mu )|\geq 2$ on $X$. 
For $\mu\in\mathcal{P}(X)$, we define $0\leq \delta (\mu )\leq 1$ to be the infimum of  
\begin{equation*}
\frac{\|\int_{X}\phi (p)\mu (dp)\|^2}{\int_{X}\|\phi (p)\|^2 \mu (dp)} 
\end{equation*}
over all mappings $\phi :\mathrm{Supp}(\mu )\to\mathcal{H}$ 
to a Hilbert space $\mathcal{H}$ 
such that 
\begin{align}
\|\phi (p)\| &=d\left( p,\mathrm{bar}(\mu )\right),\label{umbrellaborn}\\
\|\phi (p)-\phi (q)\| &\leq d(p ,q)\label{1lip}
\end{align}
for all $p,q \in\mathrm{Supp}(\mu )$. 
We define the {\em Izeki-Nayatani invariant} $\delta (X)$ of $X$ by
\begin{equation*}
\delta (X)=\sup_{\mu\in\mathcal{P} (X)}\delta (\mu ) \in\lbrack 0,1\rbrack.
\end{equation*}
\end{definition}

\begin{remark}
Notice that a mapping $\phi$ of $\mu\in\mathcal{P}(X)$ which satisfies 
\eqref{umbrellaborn} and \eqref{1lip} always exists.
To see that, fix a unit vector $e \in \mathcal{H}$.
Define $\phi(p)= d(p,\mathrm{bar}(\mu))e$.
Then by the triangle inequality, \eqref{1lip} is satisfied. 
\end{remark}

Izeki-Naytani invariant is designed to estimate the Wang invariant in comparison with the linear spectral gap. 
In \cite{IN}, Izeki and Nayatani proved the following. 

\begin{proposition}\label{delta-wang-prop}
Let $X$ be a complete $\mathrm{CAT}(0)$ space and $G$ be a weighted graph.  
Then, we have 
\begin{equation*}
\left(1-\delta (X)\right)\mu_1 (G)\leq \lambda_1 (G,X) \leq\mu_1 (G). 
\end{equation*}
\end{proposition}

It is known that there are complete $\mathrm{CAT}(0)$ spaces $X$ with $\delta (X)=1$. 
Kondo \cite{K2} constructed the first examples of such spaces.  
On the other hand, 
the present author \cite{To2} proved the following criterion for a complete $\mathrm{CAT}(0)$ space $X$ 
to be $\delta (X)<1$ (see Theorem 5.4 in \cite{To2}). 

\begin{theorem}\label{property-P-th}
Let 
$0<\theta<\frac{\pi}{2}$, $0<\alpha <1$ and $\varepsilon >0$. 
Let us say that 
a metric space $(S,d_S)$ has the \textit{property} $\mathrm{P}(\theta ,\alpha ,\varepsilon )$ if there 
exists a finite subset $S'\subset S$ 
such that 
$$
\left|\left\{ s\in S'\hspace{1mm}|\hspace{1mm} \|d_S (x,s)-d_S (y,s)\|\geq\varepsilon\right\}\right| \geq\alpha \left| S'\right|
$$
holds for 
every $x,y \in S$ with $d_S (x,y)\geq\theta$. 
Let $X$ be a complete $\mathrm{CAT}(0)$ space. 
If each tangent cone $TC_p X$ of $X$ 
is isometric to a 
(finite or infinite) product of the Euclidean cones over
metric spaces each of which has the property 
$\mathrm{P}(\theta , \alpha ,\varepsilon )$, 
then there exists a constant $C(\theta ,\alpha ,\varepsilon )<1$
depending only on 
$\theta$, $\alpha$ and $\varepsilon$ such that 
$$
\delta (X)\leq C(\theta ,\alpha ,\varepsilon ). 
$$
\end{theorem}

The following corollary is used to prove Theorem \ref{cocompact-th} in Section \ref{cocompact-sec}. 

\begin{corollary}\label{gh-precompact-corollary}
A complete $\mathrm{CAT}(0)$ space $X$ satisfies $\delta (X)<1$ 
if the family $\{ S_p X\}_{p\in X}$ consists of all spaces of directions of $X$ 
is Gromov-Hausdorff precompact.  
\end{corollary}

We recall that the Gromov-Hausdorff precompactness 
is equivalent to the uniform total boundedness 
which is defined as follows. 

\begin{definition}\label{uniformly-tatally-bounded}
A family $\mathcal{X}$ 
of metric spaces is \textit{uniformly totally bounded} 
if the following two conditions are satisfied: 
\begin{enumerate}
\item
There is a constant $D>0$ such that 
$\mathrm{diam}(X)\leq D$ for all $X\in\mathcal{X}$. 
\item
For any $\varepsilon >0$ there exists $N(\varepsilon )\in\mathbb{N}$ such that 
each $X\in\mathcal{X}$ contains a subset $S_{X,\varepsilon}\subset X$ with the following 
property: 
the cardinality of $S_{X,\varepsilon}$ is at most $N(\varepsilon )$ and 
$X$ is covered by the union of 
all open $\varepsilon$-balls whose centers are in $S_{X,\varepsilon}$. 
\end{enumerate}
\end{definition}

\begin{proof}[Proof of Corollary \ref{gh-precompact-corollary}]
It suffices to show that 
if $\mathcal{X}$ is a Gromov-Hausdorff precompact family of metric spaces, 
then there exist constants 
$0<\theta<\frac{\pi}{2}$, $0<\alpha <1$ and $\varepsilon >0$ such that 
every $X\in\mathcal{X}$ satisfies the property $\mathrm{P}(\theta ,\alpha ,\varepsilon )$. 
Since Gromov-Hausdorff precompactness is equivalent to uniform total boundedness, 
there exists an $N>0$ such that 
each $X\in\mathcal{X}$ contains a subset $S_X \subset X$ with the following 
property: 
the cardinality of $S_{X}$ is no greater than $N$ and 
$X$ is covered by the union of 
all open $\frac{\pi}{12}$-balls whose centers are in $S_{X}$. 

By the definition of the subset $S_X$, 
for any $x,y\in X$ with $d_X (x,y)\geq\frac{\pi}{3}$, 
there exist 
$s_0 ,s_1\in S_X$ such that 
\begin{align*}
&d_X (s_0, x)\geq\frac{\pi}{4},\quad d_X (s_{0},y)\leq\frac{\pi}{12},\\
&d_X (s_{1},y)\geq\frac{\pi}{4},\quad d_X (s_{1},x)\leq\frac{\pi}{12} .
\end{align*}
Hence, there exist two distinct elements $s_0 ,s_1 \in S$ such that  
\begin{align*}
\|d_{X}(x,s_0 )-d_{X}(y,s_0 )\|&\geq\frac{\pi}{6}, \\
\|d_{X}(x,s_1 )-d_{X}(y,s_1 )\|&\geq\frac{\pi}{6}, 
\end{align*}
for any $x,y\in X$ with $d_X (x,y)\geq\frac{\pi}{3}$. 
Thus each $X\in\mathcal{X}$ has the property $\mathrm{P}(\frac{\pi}{3},\frac{2}{N},\frac{\pi}{6})$. 
\end{proof}

The following corollary is used to prove Theorem \ref{doubling-th} in Section \ref{doubling-sec}. 

\begin{corollary}\label{doubling-P-corollary}
Let $X$ be a $\mathrm{CAT}(0)$ space and $p\in X$. 
Assume that the tangent cone $TC_p X$ is doubling with doubling constant $N\in\lbrack 0,\infty )$. 
Then there exist $0<\theta<\frac{\pi}{2}$, $0<\alpha <1$ 
and $\varepsilon >0$ depending only on $N$ such that the space of directions $S_p X$ at $p$ of $X$ has the 
property $\mathrm{P}(\theta ,\alpha ,\varepsilon )$. 
\end{corollary}

\begin{proof}
We assume that $N$ is a natural number. 
Since $TC_p X$ is doubling with doubling constant $N$, 
there exist closed balls $B_1 ,B_2 ,$ $\ldots$ ,$B_{N^2}$ 
with diameter at most $\frac{1}{4}$, which 
cover the closed ball of radius $1$ centered at the origin of the cone $TC_p X$. 
Hence $S_p X$ is covered by $\{\iota^{-1}\left(B_i \right)\}$, 
where $\iota :S_p X\to TC_p X$ is the canonical inclusion. 
By the inequality \eqref{inclusion}, 
each $\iota^{-1}\left( B_i \right)$ has diameter at most $\frac{\pi}{8}$. 
Thus the lemma follows from the similar argument as in the proof of corollary \ref{gh-precompact-corollary}. 
\end{proof}

Some of the known estimates of the Izeki-Nayatani invariant are: 
\begin{itemize}
\item
We have $\delta (\mathcal{H})=0$ for a Hilbert space $\mathcal{H}$ by definition. 
\item
If $X$ is a finite or infinite dimensional Hadamard manifold
or an $\mathbb{R}$-tree, then we have $\delta (Y ) = 0$ (\cite{IN}). 
\item
Let $X_p$ be the Euclidean building $PSL(3,\mathbb{Q}_p)/PSL(3,\mathbb{Z}_p)$
for each prime number $p$.  Then, we have 
$\delta (X_p ) \geq\frac{(\sqrt{p}-1)^2}{2(p-\sqrt{p}+1)}$ (\cite{IN}). 
\item
We have $\delta (X_2 ) \leq 0.4122\ldots$ (\cite{IN}). 
\item
If $X$ is a complete $\mathrm{CAT}(0)$ cube complex, then we have $\delta (X)\leq\frac{1}{2}$ (\cite{FT}). 
\end{itemize}

\section{$\mathrm{CAT}(0)$ Spaces which admit proper cocompact group actions}
\label{cocompact-sec}

In this section, we prove the following proposition.

\begin{proposition}\label{delta-cocompact-prop}
A geodesically complete $\mathrm{CAT}(0)$ space $X$ satisfies 
$\delta (X)<1$
if it admits a proper cocompact isometric action of a group. 
\end{proposition}

Combining this proposition with Proposition \ref{product-delta} in Appendix, 
we obtain the following corollary. 

\begin{corollary}\label{product-delta-cocompact}
Let $\{ X_1 ,X_2 ,\ldots ,X_n \}$ be a finite set of 
geodesically complete $\mathrm{CAT}(0)$ spaces such that 
each $X_i$ admits a proper cocompact isometric action of a group. 
Then, there exists a constant $0\leq c<1$ such that 
any $\mathrm{CAT}(0)$ space $X$ which is isometric to a (finite or infinite) product of 
copies of spaces in $\{X_1 ,X_2 ,\ldots ,X_{n}\}$ satisfies $\delta (X)\leq c$. 
\end{corollary}

Theorem \ref{cocompact-th} follows immediately from Corollary \ref{product-delta-cocompact} and Proposition \ref{delta-wang-prop}. 
Our proof of Proposition \ref{delta-cocompact-prop} consists of two lemmas. 

\begin{lemma}\label{one}
Let $X$ be a geodesically complete $\mathrm{CAT}(0)$ space. 
If there exists a positive real number $r >0$ such that 
the family $\{ B(p,r)\}_{p\in X}$ consisting of all open $r$-balls in $X$ 
is Gromov-Hausdorff precompact, then 
the family $\{ S_p X \}_{p\in X}$ consisting of all spaces of directions is 
also Gromov-Hausdorff precompact. 
\end{lemma}

\begin{proof}
Let $p\in X$ be an arbitrary point on $X$. 
We denote the canonical inclusion of $S_p X$ into $TC_p X$ by $\iota$, and 
represent the distance functions 
of $S_p X$ and $TC_p X$ by $d_S$ and $d_T$ respectively. 

Fix some $0<r'<r$. 
By the assumption, the family $\{ B(p,r)\}_{p\in X}$ is uniformly totally bounded. 
Hence, for any $\varepsilon >0$, 
there exists a positive integer $N$ which is independent of $p$ such that 
each $B(p,r)$ is covered by $N$ open balls of radius $2r' \varepsilon /\pi$. 
Then the metric sphere 
$$
S(p,r' )
=
\{
q\in X\hspace{1mm}|\hspace{1mm}
d_X (p,q)=r'
\}
\subset
B(p,r)
$$
is 
also covered by $N$ open balls of radius $2r' \varepsilon /\pi$ in $X$.

Let $F:TC_p X\to TC_p X$ be 
the mapping associating each element of $TC_p X$ represented by $(x,t)\in S_p X\times\lbrack 0,\infty)$ to 
the element represented by 
$(x,\frac{1}{r'}t)\in S_p X\times\lbrack 0,\infty)$. 
This mapping clearly satisfies 
\begin{equation}\label{scalar-multiple}
d_T (F(v),F(w))=\frac{1}{r'}d_T (v,w)
\end{equation}
for 
all $v,w\in TC_p X$. 
Then, we have $F\circ\pi_p \left( S(p,r' )\right)\subset \iota (S_p X)$, 
where $\pi_p :X\to TC_p X$ is the $1$-Lipschitz mapping defined in Definition \ref{tangent-cone}.  
By \eqref{scalar-multiple}, 
$F\circ\pi_p \left( S(p,r' )\right)$ 
can be covered by $N$ open balls of radius $2\varepsilon /\pi$ in $TC_p X$, 

Since 
each geodesic starting from $p$ can be extended up to $S(p,r')$ by geodesic completeness of $X$, 
$F\circ\pi_p \left( S(p,r')\right)$ is no other than 
$\iota\left( \left( S_p X\right)^{\circ}\right)$, and 
$F\circ\pi_p \left( S(p,r')\right)$ is dense in $\iota (S_p X)$. 
Hence, $\iota (S_p X)$ is covered by $N$ open balls of radius $2\varepsilon /\pi$ in $TC_p Y$. 
Let us denote these balls by $B_1 ,B_2 ,\ldots ,B_N$. 
Then $\{\iota^{-1}(B_i )\}_{i=1}^N$ covers $S_p X$. 
By \eqref{inclusion}, each $\iota^{-1}(B_i )$ 
is covered by an open ball in $S_p X$ of radius $\varepsilon$. 
Hence, $S_p X$ is covered by $N$ balls of radius $\varepsilon$. 
Since $\varepsilon >0$ is arbitrary, we have proved that 
$\{S_p X\}_{p\in X}$ is uniformly totally bounded. 
Thus, it is Gromov-Hausdorff precompact.  
\end{proof}

\begin{lemma}\label{two}
Let $X$ be a metric space. 
Assume that a group $\Gamma$ acts on $X$ properly and cocompactly by isometries. 
Then 
there exists some positive real number $r >0$ such that 
the family $\{B(p,r)\}_{p\in X}$
consisting of all open $r$-balls in $X$ is a Gromov-Hausdorff precompact family of metric spaces.  
\end{lemma}

\begin{proof}
Since $\Gamma$ acts on $X$ cocompactly, there exists 
a compact subset $K\subset X$ such that $\cup_{\gamma\in\Gamma}\gamma K=X$. 
Since $\Gamma$ acts on $X$ properly, 
for every $p\in K$, there exists $r_p >0$ such that the set 
$\{\gamma\in\Gamma\hspace{1mm}|\hspace{1mm}\gamma B(p,2r_p ) \cap B(p,2r_p )\neq\phi\}$ is finite. 
Let $\{ B(p_i ,r_i ) \}_{i=1}^N$ be one of finite subcovers of 
the open cover $\{ B(p,r_p )\}_{p\in K}$ of $K$. 

Although it is a well-known fact, 
we first confirm that $X$ is locally compact in this case. 
Let $q\in X$ be an arbitrary point, and let $r_0 =\min\{ r_1 ,r_2 ,\ldots ,r_N\}$. 
Observe that if there are infinitely many elements $\gamma\in\Gamma$ with 
$B(q,r_0 )\cap\gamma K\neq\emptyset$, then there exists 
some $i\in\{1,\ldots ,N\}$ with infinitely many elements $\gamma'\in\Gamma$ satisfying 
\begin{equation}\label{ballepsiloncap}
B(q,r_0 )\cap\gamma' B(p_i ,r_i )\neq\emptyset .
\end{equation}
Also, observe that 
if we can take $\gamma_1 \in\Gamma$ and $\gamma_2 \in\Gamma$ as $\gamma'$ in \eqref{ballepsiloncap}, then 
the element $\gamma_0 =\gamma_2^{-1}\gamma_1$ satisfies 
\begin{equation}\label{balls}
B(p_i ,2r_i )\cap \gamma_0 B(p_i ,2r_i )\neq\emptyset
\end{equation}
since both balls $B(x_i ,2r_i )$ and $\gamma_0 B(x_i ,2r_i )$ contain the point $\gamma_2^{-1}q$. 
Thus, if there were infinite elements $\gamma\in\Gamma$ with 
$B(q,r_0 )\cap\gamma K\neq\emptyset$, there would be infinite $\gamma_0 \in\Gamma$ 
with \eqref{balls}. 
It contradicts the definition of $r_i$. 
Thus, 
there are only finite elements $\gamma\in\Gamma$ with $B(q,r_0 )\cap\gamma K\neq\emptyset$. 
Let 
$\gamma_1 ,\ldots ,\gamma_M$ 
be all such elements. 
Then, we have 
$$
B(q ,r_0 )\subset \cup_{j=1}^M \gamma_j K. 
$$
by the definition of $K$. 
Since the right-hand side is compact, any closed ball centered at $q$ with a radius 
less than $r_0$ is compact. Hence $Y$ is locally compact.

Therefore, there exists a precompact open ball $B_p \subset X$ centered at $p$ for any $p\in K$. 
Let $\{B_i \}_i$ be a finite subcover of the open cover $\{B_p \}_{p\in K}$ of $K$, and define $U=\cup_i B_i$. 
Then, $U$ is a precompact open subset containing $K$. 

For each point $p\in K$, we define $f(p)>0$ to be 
$
f(p)
=
\sup\{\alpha >0\hspace{1mm}|\hspace{1mm}B(p,\alpha )\subset U\}
$. 
Let $q\in K$ be an arbitrary point, and let 
$\eta >0$ be an arbitrary positive real number. 
Set $\kappa =\min\{ f(q),\eta\}$. 
Then for any $q' \in B(q,\kappa )$, we have $f(q')\ge f(q)-\eta$. 
Hence, $f$ is lower semi-continuous on $K$, 
and there exists $p_0 \in K$ on which $f$ attains the minimum value of $f$.  
Set $r=f(p_0 )$.  
Then, we have $B(p,r )\subset U$ for all $p\in K$. 

Now, we show that the family $\{B(p,r)\}_{p\in X}$ of all open $r$-balls in $X$ is uniformly totally bounded. 
Let $p\in X$ be an arbitrary point, and let $\gamma\in\Gamma$ be an element which satisfies $p\in\gamma K$. 
Then, since $\gamma^{-1}B(p,r)=B(\gamma^{-1}p,r)$ is covered by $U$, 
$B(p,r)$ is covered by precompact subset $\gamma U\subset X$ which is isometric to $U$. 
Uniformly total boundedness of the family $\{B(p,r)\}_{p\in X}$ follows straightforward from this, 
which proves the lemma. 
\end{proof}

\begin{proof}[Proof of Proposition \ref{delta-cocompact-prop}]
By Lemma \ref{one} and Lemma \ref{two}, the family $\{ S_p X\}_{p\in X}$ consisting of all spaces of directions of 
geodesically complete $\mathrm{CAT}(0)$ space $X$ is Gromov-Hausdorff precompact if 
it admits a proper cocompact isometric action of a group. 
Hence, the proposition follows from Corollary \ref{gh-precompact-corollary}. 
\end{proof}

\begin{remark}
We remark that the geodesical completeness is essential in Proposition \ref{delta-cocompact-prop}. 
In \cite{K2}, Kondo constructed a sequence of locally compact $\mathrm{CAT}(0)$ cones $T_1 ,T_2 ,T_3 ,\ldots$ with 
$\lim_{i\to\infty}\delta (T_i )=1$. 
For each $i=1,2,\ldots$, let $T'_i \subset T_i$ be a closed ball of radius $\frac{1}{i}$ centered at the origin. 
Gluing $T'_1 ,T'_2 ,\ldots$ by identifying the origin of every $T'_i$, 
then the resulting space $T'$ is a compact $\mathrm{CAT}(0)$ space satisfying $\delta (T')=1$ 
although it is not geodesically complete. 
\end{remark}

\section{Ultralimits and Doubling $\mathrm{CAT}(0)$ spaces}\label{doubling-sec}

In this section, we prove the following proposition. 

\begin{proposition}\label{delta-doubling-prop}
If a complete $\mathrm{CAT}(0)$ space $X$ is uniformly locally doubling with doubling constant $N$, 
then there exists a constant $0\leq C_N <1$ depending only on $N$ which satisfies $\delta (X)<C_N$. 
\end{proposition}

Combining this proposition with Proposition \ref{product-delta} in Appendix, 
we obtain the following corollary. 

\begin{corollary}\label{product-delta-doubling}
If a complete $\mathrm{CAT}(0)$ space $X$ is isometric to 
a (finite or infinite) product 
of uniformly locally doubling $\mathrm{CAT}(0)$ spaces with a 
common doubling constant $N\in\lbrack 1,\infty )$, 
then there exists a constant $c<1$ depending only on $N$ such that $\delta (X)\leq c$.  
\end{corollary}

Theorem \ref{doubling-th} follows immediately from Corollary \ref{product-delta-cocompact} and Proposition \ref{delta-wang-prop}. 
To prove Proposition \ref{delta-doubling-prop}, 
we show that the ultralimit of a sequence of doubling length spaces with a common doubling constant 
is also doubling with 
the same doubling constant. 
First, we recall the definitions of ultrafilters and ultralimits. 
Let $I$ be a set. 
A collection $\omega\subset 2^{I}$ of subsets of $I$ is called a {\em filter} on $I$ 
if it satisfies the following conditions: 
\begin{enumerate}
\item[(a)]
$\emptyset\not\in\omega$. 
\item[(b)]
$A\in\omega$, $A\subset B$ $\Rightarrow$ $B\in\omega$. 
\item[(c)]
$A,B\in\omega$ $\Rightarrow$ $A\cap B\in\omega$. 
\end{enumerate}
An {\em ultrafilter} is a maximal filter. 
The maximality condition can be rephrased as the following condition: 
\begin{enumerate}
\item[(d)]
For any decomposition $I=A_1 \cup\cdots\cup A_m$ of $I$ into finitely many disjoint subsets $A_1 ,\ldots ,A_m$, 
$\omega$ contains exactly one of $A_1 ,\ldots ,A_m$. 
\end{enumerate}
An ultrafilter $\omega$ on $I$ is called \textit{nonprincipal} if it satisfies the following condition: 
\begin{enumerate}
\item[(e)]
For any finite subset $F\subset I, F\not\in\omega$. 
\end{enumerate}
Zorn's lemma guarantees the existence of nonprincipal ultrafilters on any infinite set $I$. 

Fix a set $I$, and an ultrafilter $\omega$ on $I$. 
For a topological space $X$, a point $p\in X$, and a mapping $f:I\to X$, we write 
\begin{equation}\label{point-limit-eq}
\omega\textrm{-}\lim_i f(i)=p
\end{equation}
if for every neighborhood $U$ of $p$, the preimage $f^{-1}(U)$ belongs to $\omega$. 
Whenever $X$ is compact and Hausdorff, for every mapping $f:I\to X$, 
there exists a unique $p\in X$ which satisfies \eqref{point-limit-eq}.

\begin{lemma}\label{ultralimit-restriction}
Fix a set $I$, an ultrafilter $\omega$ on $I$, and a subset $J\in\omega$. 
Let $X$ be a topological space, $f:I\to X$ be a mapping, and $p\in X$.  
Then, the set 
$$\omega_J =\{ K\in\omega\hspace{1mm}|\hspace{1mm}K\subset J\}$$ 
becomes an ultrafilter on $J$. 
Moreover, if we have 
$\omega_J\textrm{-}\lim_j f|_J (j)=p$ 
for the restriction $f|_J$ of $f$ to $J$, 
then we also have 
$\omega\textrm{-}\lim_i f(i)=p$
\end{lemma}

\begin{proof}
Since it is straightforward to see that $\omega_J$ is an ultrafilter on $J$, 
we only show the ``moreover" part. 
Assume that $\omega_J\textrm{-}\lim_j f|_J (j)=p$ 
holds. 
Let $U\subset X$ be an arbitrary neighborhood of $p$. 
Then by the assumption, 
$
f|_J^{-1}(U)\in\omega_J
$. 
Then $f|_J^{-1}(U)\in\omega$ by the definition of $\omega_J$. 
Hence, we have $f^{-1}(U)\in\omega$ since we have $f|_J^{-1}(U)\subset f^{-1}(U)$, 
which shows that 
$
\omega\textrm{-}\lim_i f(i)=p
$. 
\end{proof}

Fix a set $I$ and an ultrafilter $\omega$ on $I$. 
Let $\{(X_i ,d_i )\}_{i\in I}$ be a sequence of metric spaces indexed by $I$, and let 
$\prod_{i\in I} X_i$ be the set of all 
sequences $\{ p_i \}_{i\in I}$ with $p_i \in X_i$ for every $i\in I$. 
We define a relation 
$\sim$ on $\prod_{i\in I} X_i$ by declaring 
$\{ p_i \}\sim\{q_i \}$ if and only if 
$
\omega\textrm{-}\lim_i d_i (p_i ,q_i) =0, 
$
which becomes an equivalence relation. 
We denote the set of all equivalence classes by $\omega\textrm{-}\lim_i (X_i ,d_i )$, 
or simply $\omega\textrm{-}\lim_i X_i$. 
We denote the equivalence class 
represented by a sequence $\{p_i \}\in\prod_{i\in I} X_i$ by $\omega\textrm{-}\lim_i p_i$. 
We define 
the distance $d_{\omega}(p,q)$ between 
$p=\omega\textrm{-}\lim_i p_i$ and $q=\omega\textrm{-}\lim_i q_i$ by 
$$
d_{\omega}(p,q)=\omega\textrm{-}\lim_i d_i (p_i ,q_i ) \hspace{1mm}\in\lbrack 0,\infty\rbrack .
$$
Then, $(\omega\textrm{-}\lim_i (X_i ,d_i ),d_{\omega})$ becomes a 
metric space whose distance function possibly takes the value $\infty$.

\begin{definition}
Let $\omega$ be an ultrafilter on a set $I$. 
Let $\{(X_i ,d_i )\}_{i\in I}$ be a sequence of metric spaces indexed by $I$. 
We call the metric space $\left(\omega\textrm{-}\lim_i (X_i ,d_i ),d_{\omega}\right)$ defined above 
the {\em ultralimit} of $\{ (X_i ,d_i ) \}_{i\in I}$ with respect to $\omega$. 
\end{definition}

An ultralimit 
$(\omega\textrm{-}\lim_i (X_i ,d_i ),d_{\omega})$ decomposes into 
components consisting of points of mutually finite distance. 
If we are given a basepoint $p_i$ of every $X_i$, 
we can pick out the component 
consisting of points which have finite distance 
from $\omega\textrm{-}\lim_i p_i$. 
This component is a usual metric space where the distance between every pair of points is finite, and 
we denote it by $\omega\textrm{-}\lim_i (X_i ,d_i ,p_i )$.

For a sequence $\{A_i \}_{i\in I}$ of subsets $A_i\subset X_i$, 
we denote by $\omega\textrm{-}\lim_i A_i$ the subset of $\omega\textrm{-}\lim_i (X_i ,d_i )$ 
consisting of all points which are represented by sequences in $\prod_{i\in I}A_i$.

\begin{lemma}\label{ichi}
Fix a set $I$, an ultrafilter $\omega$ on $I$, and a sequence 
$\{ (X_i ,d_i )\}_{i\in I}$ of metric spaces. 
Let $\{ A_i^{(1)}\}_{i\in I}, \ldots ,\{ A_i^{(m)}\}_{i\in I}$ be sequences of subsets such that 
$A_i^{(k)}\subset X_i$ for every $k=1,\ldots ,m$ and every $i\in I$. 
Then, we have 
\begin{equation}\label{ichi-eq}
\omega\textrm{-}\lim_i \left(\bigcup_{k=1}^m A_i^{(k)}\right)
=
\bigcup_{k=1}^m\omega\textrm{-}\lim_i A_i^{(k)}. 
\end{equation}
\end{lemma}

\begin{proof}
The right-hand side of \eqref{ichi-eq} is contained in 
the left-hand side trivially. 
Let $p$ be an arbitrary point in $\omega\textrm{-}\lim_i (\cup_{k=1}^m A_i^{(k)})$ represented by 
$\{ p_i \}\in\prod_{i\in I}(\cup_{k=1}^m A_i^{(k)})$. 
For every $k\in\{ 1,2,\ldots ,m\}$, we set 
$$
I_k =\left\{ i\in I\hspace{1mm}:\hspace{1mm}
k=\min\{l\hspace{0.5mm}:\hspace{0.5mm}p_i \in A_i^{(l)}\}\right\} .
$$ 
Then, $I=I_1 \cup\cdots\cup I_m$ is a decomposition of $I$ into disjoint subsets, and 
the ultrafilter $\omega$ contains exactly one of these subsets. 
Suppose that $l\in\{ 1,2,\ldots ,m\}$ satisfies $I_l \in\omega$. 
Choose a sequence $\{q_i \}\in\prod_{i\in I}A_i^{(l)}$ 
such that $q_i =p_i$ whenever $i\in I_l$. 
Then, we have 
$$
\omega\textrm{-}\lim_i d_i (p_i ,q_i )=
\omega_{I_l}\textrm{-}\lim_i d_i (p_i ,q_i) =0
$$
by Lemma \ref{ultralimit-restriction}. 
Hence, such a sequence $\{ q_i\}\in\prod_{i\in I}A_i^{(l)}$ also represents $p$. 
Thus, we have $p\in\omega\textrm{-}\lim_i A_i^{(l)}$, which proves the lemma.  
\end{proof}

\begin{lemma}\label{ni}
Fix a set $I$, an ultrafilter $\omega$ on $I$, and a sequence 
$\{ (X_i ,d_i )\}_{i\in I}$ of length spaces. 
Let $p=\omega\textrm{-}\lim_i p_i$ be a point on the ultralimit $\omega\textrm{-}\lim_i (X_i ,d_i )$ represented by a 
sequence $\{p_i \} \in\prod_{i\in I}X_i$. 
Then, we have 
\begin{equation}\label{ultralimit-balls-eq}
\overline{B}(p,r)=\omega\textrm{-}\lim_i \overline{B}(p_{i},r) 
\end{equation}
for any $r>0$, where $\overline{B}(p,r)$ denotes the closed ball in $\omega\textrm{-}\lim (X_i ,d_i )$ of radius $r$ centered at $p$, and $\overline{B}(p_{i},r) $ denotes the closed ball in $X_i$ of radius $r$ centered at $p_i$ for each $i$. 
\end{lemma}

\begin{proof}
The right-hand side of \eqref{ultralimit-balls-eq} is contained in 
the left-hand side trivially. 
Let $q$ be an arbitrary point in the ball $\overline{B}(p,r)\subset\omega\textrm{-}\lim_i (X_i ,d_i )$ and 
let $\{q_i \}$ be a sequence representing $q$. 
We define a new sequence $\{ q'_i \}$ as follows. 
If $i\in I$ satisfies $d_i (p_i ,q_i )\le r$, we define $q'_i =q_i$. 
If $i\in I$ satisfies $d_i (p_i ,q_i )> r+1$, we define $q'_i =p_i$. 
If $i\in I$ satisfies 
$$r+\frac{1}{m+1}< d_i (p_i ,q_i )\le r+\frac{1}{m}$$
for a positive integer $m$, 
we take an arc-length parametrized path $\gamma :\lbrack 0,L\rbrack\to X_i$ of length 
$L\le r+2/m$ joining $p_i$ to $q_i$, and 
define $q'_i$ to be the point $\gamma (L-2/m)$. 
In this case, 
we have 
$
d_i (p_i ,q'_i)\le r
$ 
and 
$
d_i (q_i ,q'_i)\le\ 2/m
$. 

To prove that $q$ is contained in the right-hand side of \eqref{ultralimit-balls-eq}, 
it suffices to show that the sequence $\{ q'_i \}$ defined above satisfies 
\begin{equation}\label{same-point}
\omega\textrm{-}\lim_i d_i (q_i ,q'_i )
=
0. 
\end{equation}
Let $U\subset\mathbb{R}$ be an arbitrary neighborhood of $0\in\mathbb{R}$. 
Choose a positive integer $m$ large enough to satisfy $\overline{B}(0,\frac{2}{m})\subset U$. 
Define the subset $I_m \subset I$ by 
$$
I_m=\left\{ i\in I\hspace{1mm}\bigg|\hspace{1mm}
d_i (p_i ,q_i )\le r+\frac{1}{m}\right\}.
$$
Then, we have $I_m\in\omega$ 
since $\omega\textrm{-}\lim_i d_i (p_i ,q_i )\le r$. 
On the other hand, 
by the definition of $q'_i$, we have 
$d_i (q_i ,q'_i )\in \overline{B}(0,\frac{2}{m})$ whenever $i\in I_m$. 
Thus, 
$$
I_m \subset
\left\{ i\in I\hspace{1mm}\bigg|\hspace{1mm}
d_i (q_i ,q'_i )\in\overline{B}\left(0,\frac{2}{m}\right)\right\}
\subset
\left\{ i\in I\hspace{1mm}|\hspace{1mm}
d_i (q_i ,q'_i )\in U\right\} .
$$
Hence, 
$\left\{ i\in I\hspace{1mm}|\hspace{1mm}
d_i (q_i ,q'_i )\in U\right\}\in\omega$, 
which proves \eqref{same-point}. 
\end{proof}

We obtain the following proposition from Lemma \ref{ichi} and Lemma \ref{ni}. 

\begin{proposition}\label{doubling-ultralimit}
Fix a set $I$, 
an ultrafilter $\omega$ on $I$. 
Let $\{(X_i,d_i )\}_{i\in I}$ be a sequence of length spaces. 
If every $(X_i ,d_i)$ is doubling with a common doubling constant for every $i\in I$, 
then the ultralimit $\omega\textrm{-}\lim_i (X_i ,d_i )$ is also doubling with the same constant. 
\end{proposition}

\begin{proof}
By the assumption, there exists $N\in\mathbb{N}$ such that 
every $(X_i ,d_i )$ is doubling with doubling constant $N$. 
Fix a point $p=\omega\textrm{-}\lim_i p_i\in\omega\textrm{-}\lim_i (X_i ,d_i )$ 
and $r>0$. 
Then, 
for each $i\in I$, there exist $N$ points $p_i^{(1)},\ldots ,p_i^{(N)}\in X_i$ such that 
$$
\overline{B}(p_i ,r)\subset\bigcup_{k=1}^N\overline{B}\left(p_i^{(k)} ,\frac{r}{2}\right) , 
$$
which implies that 
\begin{equation}\label{ball-limit-eq}
\omega\textrm{-}\lim_i \overline{B}(p_i ,r)
\subset
\omega\textrm{-}\lim_i \left\{\bigcup_{k=1}^N\overline{B}\left( p_i^{(k)},\frac{r}{2}\right)\right\} .
\end{equation}
The left-hand side of \eqref{ball-limit-eq} equals $\overline{B}(p,r)$ 
by Lemma \ref{ni}, and the right-hand side equals 
$\cup_{k=1}^N\overline{B} (\omega\textrm{-}\lim_i p_i^{(k)}, r/2)$
by Lemma \ref{ichi} and \ref{ni}. 
Hence, we obtain   
$$
\overline{B}(p,r)
\subset
\bigcup_{k=1}^N\overline{B} \left(\omega\textrm{-}\lim_i p_i^{(k)}, \frac{r}{2}\right), 
$$
which proves the proposition. 
\end{proof}

\begin{proposition}\label{omega-tangent-cone}
Fix a $\mathrm{CAT}(0)$ space $(X, d_X )$, $p\in X$ and a nonprincipal ultrafilter $\omega$ on $\mathbb{N}$.   
For every $n\in\mathbb{N}$, we define another metric $d_n$ on $X$ by  
$$
d_n (p,q)= n d_X (p,q) , \quad p,q\in X. 
$$
Then the tangent cone $TC_p X$ isometrically embeds into $\omega\textrm{-}\lim_{n}(X,d_n ,p)$.
\end{proposition}

\begin{proof}
We construct an embedding 
$f:\mathrm{Cone}\left( (S_p X)^{\circ}\right)\to\omega\textrm{-}\lim_{n}(X, d_n ,p)$. 
For the origin $o\in\mathrm{Cone}\left( (S_p Y)^{\circ}\right)$, we define 
$f(o)=\omega\textrm{-}\lim_n p$. 
For $v\in\mathrm{Cone}\left( (S_p Y)^{\circ}\right)\backslash\{ o\}$, we define $f(v)$ as follows. 
Suppose that $v$ is represented by $(\lbrack\gamma\rbrack ,r)\in (S_p Y)^{\circ}\times (0,\infty )$, 
where $\lbrack\gamma\rbrack$ denotes the direction represented by a nontrivial geodesic 
$\gamma :\lbrack 0,a\rbrack\to X$ starting from $p$. 
We define a sequence $\{ p_n \}\in\prod_{n\in\mathbb{N}}X_n$ by 
$$
p_n
=
\begin{cases}
\gamma (\frac{r}{n}),\quad &\textrm{if}\hspace{3mm}\frac{r}{n}\leq a, \\
p,\quad &\textrm{if}\hspace{3mm}a<\frac{r}{n}, 
\end{cases}
$$
and define $f(v)\in\omega\textrm{-}\lim_{n}(X,d_n ,p)$ by 
$$
f(v)=\omega\textrm{-}\lim_{n}p_n .
$$
Then, by \eqref{cosine} and the definition of the distance functions on Euclidean cones, 
it is easily seen that the above definition of $f$ is independent of the choices of $\gamma$, and 
that $f$ becomes an isometry. 
Since $TC_p X$ is the metric completion of $\mathrm{Cone}\left( (S_p X)^{\circ}\right)$ and 
an ultralimit is always complete (see \cite[Chapter I, Lemma 5.53]{BH}), $f$ extends to the isometric embedding of $TC_p X$, 
which completes the proof. 
\end{proof}

Combining Proposition \ref{doubling-ultralimit} and Proposition \ref{omega-tangent-cone}, 
we obtain the following proposition.

\begin{proposition}\label{tangent-cone-doubling}
Fix $N\in\lbrack 1,\infty )$. 
Suppose that a $\mathrm{CAT}(0)$ space $(X,d_X )$ is uniformly locally doubling with doubling constant $N$. 
Then every tangent cone $TC_p X$ of $X$ is doubling with doubling constant $N$. 
\end{proposition}

\begin{proof}
Choose a nonprincipal ultrafilter $\omega$ on $\mathbb{N}$. 
Fix $p\in X$. 
For each $n\in\mathbb{N}$, let $d_n$ be a metric on $X$ defined by 
$$
d_n (p,q)= n d_X (p,q) , \quad p,q\in X. 
$$
Since $(X, d_X )$ is locally doubling with doubling constant $N$, there exists $\varepsilon >0$ such that 
the closed $\varepsilon$-ball in $(X,d_X)$ centered at $p$ is doubling with doubling constant $N$. 
Hence, for every $n$, the closed $n\varepsilon$-ball of $(X,d_n )$ centered at $p$ is doubling with 
doubling constant $N$. 
Fix an arbitrary $q=\omega\textrm{-}\lim_{n}q_n \in\omega\textrm{-}\lim_{n}(X, d_n ,p)$ and $r>0$. 
Let $s\geq 0$ be the distance between $q$ and the basepoint $\omega\textrm{-}\lim_n p$. 
We can assume that the sequence $\{q_n\}$ satisfies 
$d_n (p, q_n )\leq 2s$ for every $n$ 
by replacing all points $q_n$ with $d (p, q_n )> 2s$ by $p$ if necessary. 
Then,  
the closed $r$-ball in $(X ,d_n )$ centered at $q_n$ is 
doubling with doubling constant $N$ whenever $n\geq\frac{r+2s}{\varepsilon}$ since it is 
contained in the closed 
$n\varepsilon$-ball in $(X,d_n )$ centered at $p$. 
Hence, by Lemma \ref{ni} and Proposition \ref{doubling-ultralimit},
the closed $r$-ball in the ultralimit $\omega\textrm{-}\lim_{n}(X,d_n ,p)$ 
centered at the basepoint $\omega\textrm{-}\lim_n p$ 
is also doubling with doubling constant $N$. 
Since 
$TC_p Y$ embeds isometrically into $\omega\textrm{-}\lim_{n}(Y,d_n ,p)$ by Lemma \ref{omega-tangent-cone}, 
the proposition follows. 
\end{proof}

\begin{proof}[Proof of Proposition \ref{delta-doubling-prop}]
By Proposition \ref{doubling-ultralimit}, every tangent cone $TC_p X$ of 
uniformly locally doubling $\mathrm{CAT}(0)$ space $X$ with doubling constant $N$ is doubling with constant $N$. 
Hence, the proposition follows from Corollary \ref{gh-precompact-corollary}. 
\end{proof}

\section{Applications}\label{applications-sec}

As we described in Section \ref{intro-sec}, 
Theorem \ref{cocompact-th} and \ref{doubling-th} 
yield the fixed point property of a random group for 
a space which satisfies the hypothesis of either theorem, 
and coarse non-embeddability of sequences of expanders into such a space. 
In this section, we recall some related definitions 
and state these conclusions precisely. 

First we recall the definition of a random group of the graph model which was introduced by Gromov \cite{G}. 
A \textit{path} on a graph $G=(V,E)$ is a finite sequence 
$$
(u_1 ,u_2 ), (u_2 ,u_3),\ldots ,(u_{n-1},u_n ), (u_n ,u_{n+1} )
$$
of successive directed edges in $\vec{E}$. 
If all vertices $u_1 ,\ldots ,u_{n+1}$ are distinct, we call it an {\em embedded path}. 
A {\em cycle} is a path which starts and ends with a same vertex. 
The {\em girth} of a graph $G$ denoted by $\textrm{girth}(G)$ is the minimal length of a cycle all of whose vertices are distinct  
except the starting and ending ones. 
If there is no such a cycle in $G$, $\textrm{girth}(G)$ is defined to be $\infty$. 
The {\em diameter} of $G$ denoted by $\mathrm{diam}(G)$ is the supremum of the graph distance between two vertices in $G$. 

\begin{definition}[random groups of the graph model \cite{G}]\label{random-group-def}
Fix a positive integer $k$, and a sequence $\{G_{\ell}=(V_{\ell}, E_{\ell})\}_{\ell\in L}$ of graphs 
indexed by an unbounded set $L$ of positive integers. 
Let $\Gamma =F_k$ be the free group generated by $S=\{ s_1^{\pm},\ldots ,s_k^{\pm}\}$. 
An S-{\em labelling} of $G_{\ell}$ is a mapping $\alpha :\vec{E_{\ell}}\to S$ which satisfies 
$\alpha ((u,v))=\alpha ((v,u))^{-1}$ for every $(u,v)\in\overrightarrow{E_{\ell}}$. 
We denote by $\Lambda (G_{\ell} ,k)$ the set of all $S$-labellings of $G_{\ell}$. 
For every $\alpha\in\Lambda (G_{\ell},k )$ and every path $\vec{p}=(\vec{e}_1 ,\ldots ,\vec{e}_l )$ in $G_{\ell}$, 
we denote $\alpha (\vec{p})=\alpha (\vec{e}_1 )\cdots\alpha (\vec{e}_l )\in\Gamma$.  
We define $R_{\alpha}=\{\alpha (\vec{c} )\in\Gamma\hspace{1mm}|\hspace{1mm}\vec{c}\textrm{ is a cycle in }G_{\ell} \}$ 
and $\Gamma_{\alpha}=\Gamma /\overline{R_{\alpha}}$, 
where $\overline{R_{\alpha}}$ is the normal closure of $R_{\alpha}$. 

Now, for each group property $P$,  
we say that a {\em random group associated with} $\{ G_{\ell} \}$ {\em has property} $P$ if 
we have 
$$
\lim_{\ell\to\infty}
\frac{|\left\{\alpha\in\Lambda (G_{\ell} ,k)\hspace{1mm}|\hspace{1mm}\Gamma_{\alpha}\textrm{ has 
the property }P\right\}|}{|\Lambda (G_{\ell} ,k )|}
=1. 
$$
\end{definition}

In what follows, we fix a positive integer $k$ and every random group is based on 
the free group $F_k$ of rank $k$. 
In \cite{IKN2}, Izeki, Kondo and Nayatani proved the following fixed point theorem of 
a random group of the graph model. 
The following is a slight modification of Theorem 3.5 in \cite{IKN2}. 

\begin{theorem}[Izeki-Kondo-Nayatani \cite{IKN2}]\label{IKN-fpp}
Fix $C>0$, $d>0$ and $\lambda >0$. 
There exists $\beta >1$ which satisfies the following statement. 
Let $\mathcal{X}$ be the set of all complete $\mathrm{CAT}(0)$ spaces $X$ 
which satisfies 
$$
\lambda_1 (G, X)\geq C \mu_1 (G) 
$$
for every graph $G$. 
Let $\{G_{\ell}=(V_{\ell},E_{\ell})\}_{\ell\in L}$ be a sequence of graphs indexed by an unbounded set $L$ of positive integers 
which satisfies the following conditions: 
\begin{enumerate}
\item[$\mathrm{(a)}$]
$\lim_{\ell\to\infty}|V_{\ell} |=\infty$, 
\item[$\mathrm{(b)}$]
$2\leq\mathrm{deg}(u)\le d$ for all $\ell\in L$ and $u\in V_{\ell}$, 
\item[$\mathrm{(c)}$]
$\mu_1 (G_{\ell} )\geq \lambda$ for all $\ell\in L$, 
\item[$\mathrm{(d)}$]
There exists $c>0$ which satisfies 
$\mathrm{girth}(G_{\ell})\geq\ell$ and $\mathrm{diam}(G_{\ell})\leq c\cdot\ell$ for every $\ell\in L$. 
\item[$\mathrm{(e)}$]
There exists a constant $c' >0$ such that the number of embedded paths 
in $G_{\ell}$ of length less than $\frac{\ell}{2}$ is less than $c' \cdot\beta^{\ell /2}$. 
\end{enumerate}
Then, the random group associated with $\{G_{\ell} \}$ is non-elementary hyperbolic and 
any of its isometric action on any $X\in\mathcal{X}$ has a global fixed point. 
\end{theorem}

Combining Theorem \ref{cocompact-th} and  Theorem \ref{doubling-th} with Theorem \ref{IKN-fpp}, we obtain the following theorem. 

\begin{theorem}\label{fpp-conclusion-th}
Assume that a complete $\mathrm{CAT}(0)$ space $X$ satisfies the either of the following conditions. 
\begin{enumerate}
\item[$\mathrm{(i)}$]
$X$ is a (finite or infinite) product of 
copies of a finite number of spaces each of which is 
geodesically complete and 
admits a cocompact proper isometric action of a group. 
\item[$\mathrm{(ii)}$]
$X$ is a (finite or infinite) product of uniformly locally doubling $\mathrm{CAT}(0)$ spaces with a 
common doubling constant. 
\end{enumerate}
Then, the random group associated with a sequence $\{ G_{\ell} \}$ of graphs which satisfies the conditions $\mathrm{(a)}$,  
$\mathrm{(b)}$, $\mathrm{(c)}$, $\mathrm{(d)}$, $\mathrm{(e)}$ in Theorem \ref{IKN-fpp} is non-elementary hyperbolic and 
any of its isometric action on any $X\in\mathcal{X}$ has a global fixed point. 
\end{theorem}

\begin{definition}[sequences of expanders]
A sequence $\{G_{\ell} =(V_{\ell},E_{\ell} )\}_{\ell\in L}$ of graphs 
indexed by an unbounded set $L$ of positive integers is called 
a \textit{sequence of expanders} if it satisfies the following conditions: 
\begin{enumerate}
\item[$\mathrm{(1)}$]
$\lim_{\ell\to\infty}|V_{\ell} |=\infty$, 
\item[$\mathrm{(2)}$]
There exists $d>0$ which satisfies 
$\mathrm{deg}(u)\le d$ for all $\ell\in L$ and $u\in V_{\ell}$, 
\item[$\mathrm{(3)}$]
There exists $\lambda >0$ which satisfies 
$\mu_1 (G_{\ell} )\geq \lambda$ for all $\ell\in L$. 
\end{enumerate}
\end{definition}

So the graph sequence which is used to define the random group in Theorem \ref{IKN-fpp} and Theorem \ref{fpp-conclusion-th} is 
a sequence of expanders. 

\begin{definition}
Let $(X,d_X )$ and $(Y,d_Y )$ be metric spaces. 
A mapping $f:X\to Y$ is said to be a {\em coarse embedding} 
if there exist unbounded nondecreasing functions $\alpha ,\beta :\lbrack 0,\infty )\to\lbrack 0,\infty )$ which satisfy 
$$
\alpha (d_X (x,x' ))\leq d_Y (f (x),f (x' ) )\leq \beta (d_X (x,x' ))
$$
for every $x,x'\in X$. 
\end{definition}

Let $\{G_{\ell} =(V_{\ell} ,E_{\ell} )\}$ be a sequence of expanders, and let  
$d_{\ell}$ be the graph metric on $V_{\ell}$. 
Then, the sequence of expanders is said to be {\em embedded coarsely} into a metric space $(X,d_X )$ if 
there exist unbounded nondecreasing functions $\alpha ,\beta :\lbrack 0,\infty )\to\lbrack 0,\infty )$ and 
mappings $\left\{ f_{\ell} :V_{\ell} \to X\right\}_n$  
which satisfy 
$$
\alpha (d_{\ell} (u,v))\leq d_X (f_{\ell} (u),f_{\ell} (v) )\leq \beta (d_{\ell} (u,v))
$$
for every $\ell$ and every $u,v\in V_{\ell}$. 
Coarse embeddability of a sequence of expanders into a metric space $X$ is 
a well-known obstruction for $X$ to be embedded coarsely into a Hilbert space 
(see \cite{G} and \cite{Ro}). 
Combining Theorem \ref{cocompact-th} and  Theorem \ref{doubling-th} with the conclusion \textrm{(B)} in Section \ref{intro-sec}, 
we see that a space satisfying the hypothesis of either theorem does not have such an obstruction. 

\begin{theorem}\label{coarse-conclusion-th}
If a complete $\mathrm{CAT}(0)$ space $X$ satisfies the either condition $\mathrm{(i)}$ or $\mathrm{(ii)}$ 
in Theorem \ref{fpp-conclusion-th}, 
then a sequence of expanders does not embed coarsely into $X$. 
\end{theorem}

\appendix
\section{Some remarks on the Izeki-Nayatani invariant}\label{appendix-sec}

In this appendix, we collect some basic facts concerning the Izeki-Nayatani invariant, 
which are not mentioned in Section \ref{delta-sec}. 
First, the following proposition describes a basic behavior of the Izeki-Nayatani invariant under taking product, 
which is a slight generalization of Proposition 6.5 of \cite{IN} and quite similar to Lemma 4.3 of \cite{To2}.
We include its proof for the sake of completeness. 

\begin{proposition}\label{product-delta}
Let $X_1 ,X_2 ,X_3 ,\ldots$ be complete $\mathrm{CAT}(0)$ spaces. 
Let $X$ be a product of $X_1 ,X_2 ,X_3 ,\ldots$ (with respect to some basepoints). 
Then we have 
$$
\delta (X)=\sup \{\delta (X_i )\hspace{1mm}|\hspace{1mm} i=1,2,3,\ldots \} .
$$
\end{proposition} 

\begin{proof}
The inequality 
$
\delta (X)\geq\sup \{\delta (X_i )\hspace{1mm}|\hspace{1mm} i=1,2,3,\ldots \}
$ 
is obvious since every $X_i$ is isometrically embedded into $X$. 
Let $\mu =\sum_{i=1}^m t_i \mathrm{Dirac}_{p_i}$ be a finitely supported  probability measure 
on $X$ whose support contains at least two points. 
We write 
$p_i =(p_i^{(1)}, p_i^{(2)}, p_i^{(3)},\ldots )\in\prod_{n} X_n =X$ for each $i=1,\ldots ,n$. 
For each $n$, we define a probability measure $\mu_{n}$ on $X_n$ to be 
$$
\mu_{n}=\sum_{i=1}^m t_i \mathrm{Dirac}_{p_i^{(n)}} .
$$
Let $\mathrm{bar}(\mu )=(b_1 ,b_2 ,b_3 ,\ldots )$ be the barycenter of $\mu$. 
Then we have $\mathrm{bar}(\mu_n )=b_n$ for every $n$ since 
if we had $\mathrm{bar}(\mu_{n} )\neq b_{n}$ for some $n$, 
then it would follow that 
$$
\int_{X}d_X (p,b')^2 \mu (dp)
<
\int_{X}d_X (p,\mathrm{bar}(\mu ))^2 \mu (dp), 
$$
where $b'$ denotes the point on $X$ such that the $n$-th component is $b_n$ and the $i$-th component is 
$\mathrm{bar}(\mu_i )$ for every $i\neq n$. 

For each $n$, let $\phi_n :\mathrm{Supp}(\mu_n )\to\mathcal{H}_n$ be a realization of $\mu_n$ 
which satisfies 
$$
\delta (\mu_n )=\frac{\|\int_{X_n}\phi_n (p)\mu_n (dp)\|^2}{\int_{X_n}\|\phi_n (p)\|^2 \mu_n (dp)}. 
$$
Such a realization $\phi_n$ exists for every $n$ since the space of all realizations of $\mu_n$ is compact. 
Define a mapping $\phi :\mathrm{Supp}(\mu )\to\mathcal{H}_1 \oplus\mathcal{H}_2 \oplus\mathcal{H}_3 \oplus\cdots$ 
as 
$$
\phi(p_i )
=
\left( \phi_1 (p_i^{(1)}), \phi_2 (p_i^{(2)}), \phi_3 (p_i^{(3)}), \ldots \right), \quad i=1,\ldots ,m. 
$$
Then, it is easily seen that $\phi$ is a realization of $\mu$. 
And, we have 
\begin{multline*}
\delta(\mu )\leq
\frac{\|\int_X \phi (p)\mu (dp)\|^2}{\int_X \|\phi (p )\|^2 \mu (dp)}
=
\frac{\sum_{n=1}^{\infty}\|\sum_{i=1}^m t_i \phi_n (p_i^{(n)} )\|^2}{
\sum_{n=1}^{\infty}\sum_{i=1}^m t_i \|\phi_n (p_i^{(n)} )\|^2}
\\
\leq
\sup_n \frac{\|\sum_{i=1}^m t_i \phi_n (p_i^{(n)} )\|^2}{\sum_{i=1}^m t_i \|\phi_n (p_i^{(n)} )\|^2}
\leq
\sup_n \delta (\mu_n ),
\end{multline*}
which proves the desired inequality   
$\delta (X)\leq\sup \{\delta (X_i )\hspace{1mm}|\hspace{1mm} i=1,2,3,\ldots \}$. 
\end{proof}

Although the Izeki-Nayatani invariant is defined as a global invariant of the space,  
it can be estimated by the local property of the space. 
To see this, we define the following notation, which is introduced in 
\cite{IN}.
\begin{definition}[Izeki-Nayatani \cite{IN}]
Let $X$ be a complete $\mathrm{CAT}(0)$ space, 
and $p\in X$. 
We define $\delta (X, p)\in\lbrack 0,1\rbrack$ to be
$$
\delta (X, p )=
\sup\left\{\delta (\nu)\hspace{1mm}|\hspace{1mm}
\nu\in\mathcal{P}(X), \mathrm{bar}(\nu )=p 
\right\} ,
$$
where $\mathcal{P}(X)$ is the space of all finitely supported probability measures on Y 
whose supports contain at least two points. 
If no such $\nu$ exists, we define
$\delta(X,p)=-\infty$. 
\end{definition}

Although, the following proposition is basic and important,  
it seems that there is no reference containing its complete proof.  

\begin{proposition}\label{cone-delta}
Suppose that $X$ is a complete $\mathrm{CAT} (0)$ space. Then we have 
\begin{equation}
\delta (X)=
\sup\{\delta (TC_p X,\hspace{0.8mm} o)\hspace{1mm}|\hspace{1mm}
p\in X\}
=
\sup\{\delta (TC_p X)\hspace{1mm}|\hspace{1mm}
p\in X\} ,
\label{cone-delta-converse}
\end{equation}
where $o$ denotes the origin of the tangent cone $TC_p X$. 
\end{proposition}

\begin{proof}
The inequality 
$$
\delta (X)\leq
\sup\{\delta (TC_p X,\hspace{0.8mm} o)\hspace{1mm}|\hspace{1mm}
p\in X\}
$$
was proved in \cite[Lemma 6.2]{IN}, and the inequality 
$$
\sup\{\delta (TC_p X,\hspace{0.8mm} o)\hspace{1mm}|\hspace{1mm}p\in X\}
\leq
\sup\{\delta (TC_p X)\hspace{1mm}|\hspace{1mm}
p\in X\}
$$
is trivial from the definition. 
So we need only to prove the inequality 
\begin{equation}\label{T<=X}
\sup\{\delta (TC_p X)|p\in X\}\leq\delta (X).
\end{equation}
To this end, 
it suffices to show that $\delta (TC_p X)\leq\delta (X)$ for 
any $p\in X$. 
By Proposition 4.2 of \cite{IKN}, if 
$\{ X_n\}_{n\in\mathbb{N}}$ is a sequence of complete $\mathrm{CAT}(0)$ spaces, 
$\omega$ is a nonprincipal ultrafilter on $\mathbb{N}$, and 
$X_{\omega}$ is the ultralimit of $\{ X_n\}$ with respect to $\omega$, 
then $\delta (X_{\omega})\leq\omega\textrm{-}\lim_{n}\delta (X_n )$ holds. 
Combining this with Proposition \ref{omega-tangent-cone} in Section \ref{doubling-sec}, the inequality \eqref{T<=X}
follows immediately. 
\end{proof}

\begin{remark}
Although the previous version of this paper \cite{To3} also includes 
the proof of Proposition \ref{cone-delta} without using the notion of ultralimit, 
we omit it here to avoid redundancy. 
\end{remark}


\begin{thebibliography}{10}

\bibitem{BH}
\textsc{M.~R. Bridson and A.~Haefliger}, 
\newblock Metric spaces of non-positive curvature, 
volume 319 of {\em
  Grundlehren der Mathematischen Wissenschaften [Fundamental Principles of
  Mathematical Sciences]}.
\newblock Springer-Verlag, Berlin, 1999.

\bibitem{BBI}
\textsc{Burago~Y. Burago, D. and S.~Ivanov}, 
\newblock A course in metric geometry, 
volume~33 of {\em Graduate studies
  in Math.}
\newblock Amer. Math. Soc., Providence, RI, 2001.

\bibitem{FT} \textsc{K.~Fujiwara and T.~Toyoda}, 
\newblock Random groups have fixed points on CAT(0) cube complexes,  
\newblock Proc. Amer. Math. Soc., \textbf{140} (2012), 1023--1031. 

\bibitem{G} \textsc{M.~Gromov}, 
\newblock Random walk in random groups, 
\newblock Geom. Funct. Anal., \textbf{13}(1) (2003), 73--146.

\bibitem{IKN} \textsc{H.~Izeki, T.~Kondo, and S.~Nayatani}, 
\newblock Fixed-point property of random groups.
\newblock Ann. Global Anal. Geom., \textbf{35}(4) (2009), 363--379.

\bibitem{IKN2} \textsc{H.~Izeki, T.~Kondo, and S.~Nayatani}, 
\newblock $N$-step energy of maps and fixed-point property of random groups, 
\newblock Groups Geom. Dyn. \textbf{6}(4) (2012), 701--736. 

\bibitem{IN} \textsc{H.~Izeki and S.~Nayatani}, 
\newblock Combinatorial harmonic maps and discrete-group actions on {H}adamard spaces, 
\newblock {\em Geom. Dedicata}, \textbf{114} (2005), 147--188.

\bibitem{K2}
\textsc{T.~Kondo}, 
\newblock ${\rm CAT}(0)$ spaces and expanders, 
\newblock Math. Z. \textbf{271} (2012),  no. 1-2, 343--355. 

\bibitem{NS}
\textsc{A.~Naor and L.~Silberman}, 
\newblock Poincar\'e  inequalities, embeddings, and wild groups, 
\newblock Compositio Mathematica \textbf{147} (2011), no. 5, 1546--1572.

\bibitem{P} \textsc{P.~Pansu}, 
\newblock Superrigidit\'e geometrique et applications harmoniques, 
\newblock S\'eminaires et congr\`es 18, 375--422, 
Soc. Math. France, Paris, 2008. 

\bibitem{Ro} \textsc{J.~Roe}, 
\newblock Lectures on coarse geometry, 
volume~31 of {\em University
  Lecture Series}.
\newblock Amer. Math. Soc., Providence, RI, 2003.

\bibitem{St} \textsc{K.T. Sturm}, 
\newblock Probability measures on metric spaces of nonpositive curvature, 
\newblock {Heat kernels and analysis on manifolds, graphs, and metric
  spaces} (Paris, 2002), 357--390, 2003.
\newblock {\em Contemp. Math}., 338, Amer. Math. Soc., Providence, RI, 2003. 

\bibitem{To2} \textsc{T.~Toyoda}, 
\newblock ${\rm CAT}(0)$ spaces on which certain type of singularity is bounded, 
\newblock Kodai Math. J., \textbf{33} (2010), 398--415.

\bibitem{To3} \textsc{T.~Toyoda}, 
\newblock Fixed point property for a $\mathrm{CAT}(0)$ space which admits a proper cocompact group action, 
\newblock preprint, arXiv:1102.0729v2, 2011.

\bibitem{Wa} \textsc{M.-T. Wang}, 
\newblock Generalized harmonic maps and representations of discrete groups, 
\newblock {\em Comm. Anal. Geom.}, \textbf{8}(3) (2000), 545--563.

\end{thebibliography}
\end{document}